# ON OPTIMALITY OF BAYESIAN TESTIMATION IN THE NORMAL MEANS PROBLEM

By Felix Abramovich, Vadim Grinshtein and Marianna Pensky[1]

*Tel Aviv University, The Open University of Israel and
University of Central Florida*

We consider a problem of recovering a high-dimensional vector $\mu$ observed in white noise, where the unknown vector $\mu$ is assumed to be sparse. The objective of the paper is to develop a Bayesian formalism which gives rise to a family of $l_0$-type penalties. The penalties are associated with various choices of the prior distributions $\pi_n(\cdot)$ on the number of nonzero entries of $\mu$ and, hence, are easy to interpret. The resulting Bayesian estimators lead to a general thresholding rule which accommodates many of the known thresholding and model selection procedures as particular cases corresponding to specific choices of $\pi_n(\cdot)$. Furthermore, they achieve optimality in a rather general setting under very mild conditions on the prior. We also specify the class of priors $\pi_n(\cdot)$ for which the resulting estimator is adaptively optimal (in the minimax sense) for a wide range of sparse sequences and consider several examples of such priors.

**1. Introduction.** Consider a problem of estimation of a high-dimensional multivariate Gaussian mean with independent terms and common variance,

$$(1.1) \qquad y_i = \mu_i + \sigma z_i, \qquad z_i \overset{i.i.d.}{\sim} N(0,1), \qquad i=1,\ldots,n.$$

The variance $\sigma^2$ is assumed to be known and the goal is to estimate the unknown mean vector $\mu$ from a set $\Theta_n \subset \mathbb{R}^n$. This is a well-studied problem that arises in various statistical settings, for example, model selection or orthonormal regression.

Some extra assumptions are usually placed on $\Theta_n$. We assume that the vector $\mu$ is *sparse*, that is, most of its entries are zeroes or "negligible" and only a small fraction is "significantly large." The indices of the large entries

Received February 2006; revised December 2006.
[1]Supported in part by NSF Grants DMS-05-05133 and DMS-06-52524.

are, however, not known in advance. Formally, the sparsity assumption can be quantified in terms of so-called nearly-black objects (Donoho et al. [12]) or strong and weak $l_p$-balls (Johnstone [20], Donoho and Johnstone [10, 11]) discussed below. In sparse cases the natural estimation strategy is thresholding.

It is well known that various thresholding rules can be considered as penalized likelihood estimators minimizing

$$(1.2) \qquad \|y - \mu\|_2^2 + P(\mu)$$

for the corresponding penalties $P(\mu)$. The traditional $l_2$ penalty $P(\mu) = \lambda_n^2 \|\mu\|_2^2$ does not lead to a thresholding but to a linear shrinkage estimator $\hat{\mu}_i^* = \frac{1}{1+\lambda_n^2} y_i$. The $l_1$ penalty produces a "shrink-or-kill" soft thresholding, where $\hat{\mu}_i^* = sign(y_i)(|y_i| - \lambda_n/2)_+$, which coincides with the LASSO estimator of Tibshirani [25]. The general $l_p$, $p > 0$, penalty yields bridge regression (Frank and Friedman [16]) and results in thresholding when $p \leq 1$. Wider classes of penalties leading to various thresholding rules are discussed in Antoniadis and Fan [6], Fan and Li [13] and Hunter and Li [19].

All the penalties mentioned above are related to *magnitudes* of $\mu_i$. In this paper we consider the $l_0$, or complexity type penalties, where the penalty is placed on the *number of nonzero* $\mu_i$. The $l_0$ quasi-norm of a vector $\mu$ is defined as the number of its nonzero entries, that is, $\|\mu\|_0 = \#\{i : \mu_i \neq 0\}$. In the simplest case, the complexity penalty $P(\mu) = \lambda_n^2 \|\mu\|_0$ and minimization of (1.2) obviously result in minimizing

$$(1.3) \qquad \sum_{i=k+1}^{n} y_{(i)}^2 + \lambda_n^2 k$$

over $k$, where $|y|_{(1)} \geq \cdots \geq |y|_{(n)}$. Such a procedure implies a "keep-or-kill" hard thresholding with a (fixed) threshold $\lambda_n$,

$$\hat{\mu}_i^* = y_i I\{|y_i| \geq \lambda_n\}, \qquad i = 1, \ldots, n.$$

The widely-known universal threshold of Donoho and Johnstone [9] is $\lambda_U = \sigma\sqrt{2\ln n}$ and, as $n \to \infty$, the resulting estimator comes within a constant factor of asymptotic minimaxity for $l^r$ losses simultaneously throughout a range of various sparsity classes (Donoho and Johnstone [10, 11]).

A complexity penalization of type (1.3) is closely connected to model selection. For example, the Akaike's [5] AIC model selection rule takes $\lambda_n = \sqrt{2}\sigma$, the Schwarz [24] BIC criterion corresponds to $\lambda_n = \sigma\sqrt{\ln n}$, while the RIC criterion of Foster and George [14] adjusted for (1.1) implies $\lambda_n = \sigma\sqrt{2\ln n}$.

A natural extension of (1.3) is to consider a variable penalization sequence $\lambda_{i,n}$, that is,

$$(1.4) \qquad \|y - \mu\|^2 + \sum_{i=0}^{\|\mu\|_0} \lambda_{i,n}^2.$$

Let $\hat{k}$ be a minimizer of

$$\sum_{i=k+1}^{n} y_{(i)}^2 + \sum_{i=0}^{k} \lambda_{i,n}^2 \tag{1.5}$$

over $k$. The resulting minimizer $\hat{\mu}^*$ of (1.4) is obviously a hard thresholding rule with a variable threshold $\lambda_{\hat{k},n} : \hat{\mu}_i^* = y_i I\{|y_i| \geq \lambda_{\hat{k},n}\}$. If $\hat{k} = 0$, all $y_i$ are thresholded and $\hat{\mu}^* \equiv 0$.

Several variable penalty estimators of the type (1.4) have been proposed in the literature. The FDR-thresholding rule of Abramovich and Benjamini [2] corresponds to $\lambda_{i,n} = \sigma z(1-(i/n)(q/2)) \sim \sigma\sqrt{2\ln((n/i)(2/q))}$, where $z(\cdot)$ are standard Gaussian quantiles and $q$ is the tuning parameter of the FDR procedure. Abramovich et al. [4] showed that, for $q \leq 1/2$, the FDR estimator achieves *sharp* (with a right constant) asymptotic minimaxity, simultaneously over an entire range of nearly black sets and $l_p$-balls with respect to $l^r$ losses. Foster and Stine [15] suggested $\lambda_{i,n} = \sigma\sqrt{2\ln(n/i)}$ from information theory considerations. The covariance inflation criterion for model selection of Tibshirani and Knight [26] adjusted for (1.1) corresponds to $\lambda_{i,n} = 2\sigma\sqrt{\ln(n/i)}$.

For a general $l_0$-type penalty $P_n(\|\mu\|_0)$, the corresponding penalized estimator $\hat{\mu}^*$ is a hard thresholding rule with the data-dependent threshold $\hat{\lambda} = |y|_{(\hat{k})}$, where $\hat{k}$ is the minimizer of

$$\sum_{i=k+1}^{n} y_{(i)}^2 + P_n(k). \tag{1.6}$$

Obviously, (1.5) can be viewed as a particular case of (1.6).

A series of recent papers has considered the $2k\ln(n/k)$-type penalties of the form

$$P_n(k) = 2\sigma^2 \zeta k(\ln(n/k) + c_{k,n}), \tag{1.7}$$

where $\zeta > 1$ and $c_{k,n}$ is a negligible term relative to $\ln(n/k)$ for $k \ll n$ (sparse cases) (e.g., Birgé and Massart [8], Johnstone [21] and Abramovich et al. [4]).

A wide class of $l_0$-type penalties satisfying certain technical conditions was considered in Birgé and Massart [8]. However, most of their results on optimality have been obtained for a particular $2k\ln(n/k)$-type penalty, while it remains somewhat unclear how to construct "meaningful" penalties from their class in general.

The objective of this paper is to develop a framework for $l_0$-penalization which is general and meaningful at the same time. The Bayesian approach provides a natural interpretation of the penalized likelihood estimators by relating the penalty models to the corresponding prior distribution on $\mu$. For the model (1.1), the penalty term in (1.2) is then proportional to the logarithm of the prior. From a Bayesian view, minimization of (1.2) corresponds to the maximum a posteriori (MAP) rule and the resulting penalized estimator (called MAP thereafter) is the posterior mode. The $l_p$-type penalties for $p > 0$ correspond to placing priors on the *magnitudes* of $\mu_i$, while the $l_0$-type penalties necessarily involve a prior on the number of nonzero $\mu_i$.

In this paper we develop a Bayesian formalism which gives rise to a family of $l_0$-type penalties in (1.6). This family is associated with various choices of the prior distributions on the number of nonzero entries of the unknown vector and, hence, is easy to interpret. Moreover, under mild conditions, the penalties considered in this paper fall within the class considered in Birgé and Massart [8], which allows us to establish optimality of the corresponding Bayesian estimators in a rather general setting. We then demonstrate that in the case when the vector $\mu$ is sparse, the MAP estimators achieve optimal convergence rates. We also specify the class of prior distributions for which the resulting estimators are adaptive for a wide range of sparse sequences and provide examples of such priors.

The paper is organized as follows. In Section 2 we introduce the Bayesian MAP "testimation" procedure leading to penalized estimators (1.6). In Section 3 we derive upper bounds for their quadratic risk and compare it with that of an ideal oracle estimator (oracle inequality). Asymptotic optimality of the proposed MAP "testimators" in various sparse settings is established in Section 4. Several specific priors $\pi_n(k)$ are considered in Section 5 as examples. In Section 6 we present a short simulation study to demonstrate the performance of MAP estimators and compare them with several existing counterparts. Some concluding remarks are given in Section 7. All the proofs are placed in the Appendix.

## 2. MAP testimating procedure.

2.1. *Thresholding as testimation.* Abramovich and Benjamini [2] demonstrated that thresholding can be viewed as a multiple hypothesis testing procedure, where, given the data $y = (y_1, \ldots, y_n)'$ in (1.1), one first simultaneously tests $\mu_i, i = 1, \ldots, n$, for significance. Those $\mu_i$'s which are concluded to be significant are estimated by the corresponding $y_i$, while nonsignificant $\mu_i$'s are discarded. Such a *testimation* procedure obviously mimics a hard thresholding rule.

In particular, the likelihood ratio test rejects the null hypothesis $H_{0i}: \mu_i = 0$ if and only if $|y|_i > \lambda_n$, where controlling the familywise error at level $\alpha$ by the Bonferroni approach leads to $\lambda_n = \sigma z(1 - \alpha/(2n)) \sim \sigma\sqrt{2\ln n} = \lambda_U$ for any reasonable $\alpha$. In other words, universal thresholding may be viewed as a Bonferroni multiple testing procedure with familywise error level of approximately $1/\sqrt{\ln n}$, which slowly approaches zero as $n$ increases. Such a severe error control explains why universal thresholding is so conservative. A less stringent alternative to a familywise error control is the false discovery

rate (FDR) criterion of Benjamini and Hochberg [7]. The corresponding FDR thresholding was considered in Abramovich and Benjamini [2, 3] in the context of wavelet series estimation and comprehensively developed further in Abramovich et al. [4] for a general normal means problem setting.

In this paper we shall follow a more general testimation approach to thresholding based on the multiple hypothesis testing procedure introduced by Abramovich and Angelini [1], which efficiently utilizes the Bayesian framework. We shall review this approach in the following section.

2.2. *MAP multiple testing procedure.* For the model (1.1), consider the multiple hypothesis testing problem, where we wish to simultaneously test $H_{0i} : \mu_i = 0$ against $H_{1i} : \mu_i \neq 0$, $i = 1, \ldots, n$.

A configuration of true and false null hypotheses is uniquely defined by the indicator vector $x$, where $x_i = I\{\mu_i \neq 0\}, i = 1, \ldots, n$. Let $k = x_1 + \cdots + x_n = \|\mu\|_0$ be the number of significant $\mu_i$ (false nulls). Assume some prior distribution $k \sim \pi_n(k) > 0, k = 0, \ldots, n$. For a given $k$ there are $\binom{n}{k}$ various configurations of true and false null hypotheses. Assume all of them to be equally likely a priori, that is, conditionally on $k$,

$$P\left(x \,\Big|\, \sum_{i=1}^n x_i = k\right) = \binom{n}{k}^{-1}.$$

Naturally, $\mu_i | x_i = 0 \sim \delta_0$, where $\delta_0$ is a probability atom at zero. To complete the prior, assume $\mu_i | x_i = 1 \sim N(0, \tau^2)$.

For the proposed hierarchical prior, the posterior distribution of configurations is given by

$$(2.8) \qquad \pi_n(x, k|y) \propto \binom{n}{k}^{-1} \pi_n(k) I\left\{\sum_{i=1}^n x_i = k\right\} \prod_{i=1}^n (B_i^{-1})^{x_i},$$

where the Bayes factor $B_i$ of $H_{0i}$ is

$$(2.9) \qquad B_i = \sqrt{1+\gamma} \exp\left\{-\frac{y_i^2}{2\sigma^2(1+1/\gamma)}\right\}$$

and the variances ratio $\gamma = \tau^2/\sigma^2$ (Abramovich and Angelini [1]).

Given the posterior distribution $\pi_n(x, k|y)$, we apply a maximum a posteriori (MAP) rule to choose the most likely configuration of true and false nulls. Generally, to find the posterior mode of $\pi_n(x, k|y)$, one should look through all $2^n$ possible configurations. However, for the proposed model, the number of candidates for a mode is, in fact, reduced to $n+1$ only. Indeed, let $\hat{x}(k)$ be a maximizer of (2.8) for a *fixed* $k$ that indicates the most plausible configuration with $k$ false null hypotheses. From (2.8) it follows immediately that $\hat{x}_i(k) = 1$ for the $k$ tests with the smallest Bayes factors $B_i$ and is zero otherwise. Due to the monotonicity of $B_i$ in $|y|_i$ [see (2.9)], this is equivalent to $\hat{x}_i(k) = 1$ corresponding to the $k$ largest $|y|_i$ and zero for the others. The Bayesian MAP multiple testing procedure then leads to finding $\hat{k}$ that maximizes

$$\ln \pi_n(\hat{x}(k), k|y) = \sum_{i=1}^{k} y_{(i)}^2 + 2\sigma^2(1+1/\gamma) \ln\left\{ \binom{n}{k}^{-1} \pi_n(k)(1+\gamma)^{\frac{-k}{2}} \right\} + const.$$

or, equivalently, minimizes

$$\sum_{i=k+1}^{n} y_{(i)}^2 + 2\sigma^2(1+1/\gamma) \ln\left\{ \binom{n}{k} \pi_n^{-1}(k)(1+\gamma)^{k/2} \right\}.$$

The $\hat{k}$ null hypotheses corresponding to $|y|_{(1)}, \ldots, |y|_{(\hat{k})}$ are rejected. The resulting Bayesian testimation yields hard thresholding with a threshold $\hat{\lambda}_{\mathrm{MAP}} = |y|_{(\hat{k})}$:

$$\hat{\mu}_i^* = \begin{cases} y_i, & |y_i| \geq \hat{\lambda}_{\mathrm{MAP}}, \\ 0, & \text{otherwise}, \end{cases}$$

and is, in fact, the posterior mode of the joint distribution $(\mu, x, k|y)$.

From a frequentist view, the above MAP estimator $\hat{\mu}^*$ is a penalized likelihood estimator (1.6) with the complexity penalty

(2.10) $$P_n(k) = 2\sigma^2(1+1/\gamma) \ln\left\{ \binom{n}{k} \pi_n^{-1}(k)(1+\gamma)^{k/2} \right\}.$$

Rewriting $\binom{n}{k} = \prod_{i=1}^{k}(n-i+1)/i$ and $\pi_n(k) = \pi_n(0)\prod_{i=1}^{k}\pi_n(i)/\pi_n(i-1)$, (2.10) yields

$$P_n(k) = 2\sigma^2(1+1/\gamma)\left( \ln \pi_n^{-1}(0) + \sum_{i=1}^{k} \ln\left\{ \frac{n-i+1}{i} \frac{\pi_n(i-1)}{\pi_n(i)} \sqrt{1+\gamma} \right\} \right)$$
(2.11)
$$= \sum_{i=0}^{k} \lambda_{i,n},$$

where $\lambda_{0,n} = 2\sigma^2(1+1/\gamma)\ln \pi_n^{-1}(0)$, $\lambda_{i,n} = 2\sigma^2(1+1/\gamma)\ln(\frac{n-i+1}{i}\frac{\pi_n(i-1)}{\pi_n(i)} \times \sqrt{1+\gamma})$, $i = 1, \ldots, n$.

In such a form the penalty (2.11) is similar to (1.5) with the notation $\lambda_{i,n}$ instead of $\lambda_{i,n}^2$ since some of them might be negative in a general case.

A specific form of the resulting Bayesian hard thresholding rule depends on the choice of a prior $\pi_n(k)$. In particular, the binomial prior $B(n, \xi_n)$ yields a *fixed* threshold $2\sigma^2(1+1/\gamma)\ln(\frac{1-\xi_n}{\xi_n}\sqrt{1+\gamma}) \sim 2\sigma^2 \ln(\frac{1-\xi_n}{\xi_n}\sqrt{\gamma})$ for sufficiently large $\gamma$. The AIC criterion corresponds to $\xi_n \sim \sqrt{\gamma}/(e+\sqrt{\gamma})$, while $\xi_n = 1/n$ leads to the universal thresholding of Donoho and Johnstone [9]. Abramovich and Angelini [1] showed that the "reflected" truncated Poisson distribution $\pi_n(k) \propto (n - \lambda_n)^{n-k}/(n-k)!$, with $\lambda_n = o(n)$ satisfying $\lambda_n/\sqrt{n \ln n} \to \infty$, approximates the FDR thresholding procedures of Benjamini and Hochberg [7] and Sarkar [23] with the FDR parameter $q_n \sim (\sqrt{\pi \gamma \ln(\sqrt{\gamma} n/\lambda_n)})^{-1}$.

REMARK 1. There is an intriguing parallel between the penalty $P_n(k)$ in (2.10) and $2k \ln(n/k)$-type penalties (1.7) introduced above. For $k \ll n$, $\ln \binom{n}{k} \sim k \ln(n/k)$ and the penalty (2.10) is of $2k \ln(n/k)$-type, where $\zeta = (1 + 1/\gamma)$ and $c_{k,n} = (1/k) \ln \pi_n^{-1}(k) + (1/2) \ln(1+\gamma)$ is defined by the choice of the prior $\pi_n(k)$. The $2k \ln(n/k)$-type penalty can be viewed, therefore, as a particular case of the more general penalty (2.10) for $\pi_n(k)$ satisfying $c_{k,n} = O(\ln(n/k))$, or, equivalently $\ln \pi_n(k) = O(k \ln(k/n))$ for $k \ll n$. Throughout the paper we discuss the relations between $P_n(k)$ and $2k \ln(n/k)$-type penalties in more detail.

In what follows, we study optimality of the proposed thresholding MAP estimators.

**3. Oracle inequality.** In this section we derive an upper bound for the quadratic risk $\rho(\hat{\mu}^*, \mu) = E\|\hat{\mu}^* - \mu\|^2$ of the MAP thresholding estimator and compare it with the ideal risk of an oracle estimator.

ASSUMPTION (A). Assume that

$$\pi_n(k) \leq \binom{n}{k} e^{-c(\gamma)k}, \qquad k = 0, \ldots, n, \tag{3.12}$$

where $c(\gamma) = 8(\gamma + 3/4)^2$.

THEOREM 1 (Upper bound). *Under Assumption* (A),

$$\rho(\hat{\mu}^*, \mu) \leq c_0(\gamma) \inf_{0 \leq k \leq n} \left\{ \sum_{i=k+1}^n \mu_{(i)}^2 \right.$$

$$\left. + 2\sigma^2 (1 + 1/\gamma) \ln\left( \binom{n}{k} \pi_n(k)^{-1} (1+\gamma)^{k/2} \right) \right\} \tag{3.13}$$

$$+ c_1(\gamma) \sigma^2$$

*for some $c_0(\gamma)$ and $c_1(\gamma)$ depending only on $\gamma$.*

The obvious inequality $\binom{n}{k} \geq (n/k)^k$ implies that, for any $\pi_n(k)$, (3.12) holds for all $k \leq ne^{-c(\gamma)}$. Applying the arguments very similar to those in the proof of Theorem 1, one gets then a somewhat weaker general upper bound for the quadratic risk of $\hat{\mu}^*$ without the requirement of Assumption (A):

COROLLARY 1. *For any prior $\pi_n(\cdot)$,*

$$\rho(\hat{\mu}^*, \mu) \leq c_0(\gamma) \inf_{0 \leq k \leq ne^{-c(\gamma)}} \left\{ \sum_{i=k+1}^{n} \mu_{(i)}^2 + 2\sigma^2(1 + 1/\gamma) \ln\left( \binom{n}{k} \pi_n(k)^{-1}(1+\gamma)^{k/2} \right) \right\}$$
$$+ c_1(\gamma)\sigma^2,$$

*where $c(\gamma) = 8(\gamma + 3/4)^2$, and $c_0(\gamma)$ and $c_1(\gamma)$ depend only on $\gamma$.*

Note that the upper bounds in Theorem 1 and Corollary 1 are nonasymptotic and hold for any $\mu \in \mathbb{R}^n$.

In order to assess the quality of the upper bound in (3.13), we compare the quadratic risk of the MAP estimator with that of the ideal estimator $\hat{\mu}_{\text{oracle}}$ which one could obtain if one had available an oracle which reveals the true vector $\mu$. This ideal risk is known to be

$$(3.14) \qquad \rho(\hat{\mu}_{\text{oracle}}, \mu) = \sum_{i=1}^{n} \min(\mu_i^2, \sigma^2)$$

(Donoho and Johnstone [9]). The ideal estimator $\hat{\mu}_{\text{oracle}}$ is obviously unavailable but can be used as a benchmark for the risk of other estimators. Note that the risk of $\hat{\mu}_{\text{oracle}}$ is zero when $\mu \equiv 0$ and, evidently, no estimator can achieve this risk bound in this case. An additional (usually negligible) term $\sigma^2$, which is, in fact, an error of unbiased estimation of one extra parameter, is usually added to the ideal risk in (3.14) for a proper comparison (see, e.g., Donoho and Johnstone [9]).

Define

$$(3.15) \qquad \begin{aligned} L_{0,n} &= 2\ln \pi_n^{-1}(0), \\ L_{k,n} &= (1/k)\ln\left( \binom{n}{k} \pi_n^{-1}(k) \right), \qquad k = 1, \ldots, n, \end{aligned}$$

and let $L_n^* = \max_{0 \leq k \leq n} L_{k,n}$. The following theorem states that the MAP thresholding estimator performs within a factor of $c_2(\gamma)(2L_n^* + \ln(1+\gamma))$ with respect to the oracle.

THEOREM 2 (Oracle inequality). *Consider the MAP thresholding estimator $\hat{\mu}^*$ and the corresponding penalty $P_n(k)$ defined in (2.10). Under Assumption (A),*

$$\rho(\hat{\mu}^*, \mu) \leq c_2(\gamma)(2L_n^* + \ln(1+\gamma))(\rho(\hat{\mu}_{\text{oracle}}, \mu) + \sigma^2)$$

*for some $c_2(\gamma)$ depending only on $\gamma$.*

To understand how tight the factor $c_2(\gamma)(2L_n^* + \ln(1+\gamma))$ is, recall that when $n$ is large, there is a sharp upper bound for the quadratic risk,

$$(3.16) \quad \inf_{\tilde{\mu}} \sup_{\mu} \frac{\rho(\tilde{\mu}, \mu)}{\rho(\hat{\mu}_{\text{oracle}}, \mu) + \sigma^2} = 2\log n(1 + o(1)) \qquad \text{as } n \to \infty$$

(Donoho and Johnstone [9]).

Therefore, no available estimator has a risk smaller than within the factor $2\log n$ from an oracle. Obvious calculus shows that if $\pi_n(k) \geq n^{-ck}, k = 1, \ldots, n$, and $\pi_n(0) \geq n^{-c}$ for some constant $c > 0$, then $L_n^* = O(\log n)$ and the MAP estimator achieves the minimal possible risk among all available estimators (3.16) up to a constant factor depending on $\gamma$:

COROLLARY 2. *Let $\pi_n(k)$ satisfy Assumption (A) and, in addition, $\pi_n(k) \geq n^{-ck}, k = 1, \ldots, n$, and $\pi_n(0) \geq n^{-c}$ for some constant $c > 0$. The resulting MAP estimator $\hat{\mu}^*$ satisfies*

$$(3.17) \quad \sup_{\mu} \frac{\rho(\hat{\mu}^*, \mu)}{\rho(\hat{\mu}_{\text{oracle}}, \mu) + \sigma^2} = c_3(\gamma) 2 \log n(1+o(1)) \qquad \text{as } n \to \infty$$

*for some $c_3(\gamma) \geq 1$.*

In particular, Corollary 2 holds for $\ln \pi_n(k) = O(k \ln(k/n))$ corresponding to the $2k \ln(n/k)$-type penalties (see Remark 1 in Section 2.2). However, the condition $\pi_n(k) \geq n^{-ck}$ required in the corollary is much weaker and covers a far wider class of possible priors.

**4. Minimaxity and adaptivity in sparse settings.** The results of Section 3 hold for any $\mu \in \mathbb{R}^n$. In this section we show that they can be improved if an extra *sparsity* constraint on $\mu$ is added. We start by introducing several possible ways to quantify sparsity and then derive conditions on the prior $\pi_n(\cdot)$ which imply asymptotic minimaxity of the resulting MAP estimator $\hat{\mu}^*$ over various sparse settings.

4.1. *Sparsity.* The most intuitive measure of sparsity is the number of nonzero components of $\mu$, or its $l_0$ quasi-norm: $\|\mu\|_0 = \#\{i : \mu_i \neq 0, i = 1, \ldots, n\}$. Define a $l_0$-ball $l_0[\eta]$ of standardized radius $\eta$ as a set of $\mu$ with at most a proportion $\eta$ of nonzero entries, that is,

$$l_0[\eta] = \{\mu \in \mathbb{R}^n : \|\mu\|_0 \leq \eta n\}.$$

In a wider sense sparsity can be defined by the proportion of large entries. Formally, define a weak $l_p$-ball $m_p[\eta]$ with standardized radius $\eta$ as

$$m_p[\eta] = \{\mu \in \mathbb{R}^n : |\mu|_{(i)} \leq \eta(n/i)^{1/p}, i = 1, \ldots, n\}.$$

Sparsity can be also measured in terms of the $l_p$-norm of a vector. A strong $l_p$-ball $l_p[\eta]$ with standardized radius $\eta$ is defined as

$$l_p[\eta] = \left\{\mu \in \mathbb{R}^n : \frac{1}{n}\sum_{i=1}^n |\mu_i|^p \leq \eta^p\right\}.$$

There are important relationships between these balls. The $l_p$-norm approaches $l_0$ as $p$ decreases. The weak $l_p$-ball contains the corresponding strong $l_p$-ball, but only just:

$$l_p[\eta] \subset m_p[\eta] \not\subset l_{p'}[\eta], \qquad p' > p.$$

The smaller $p$ is, the sparser is $\mu$. Sparse settings correspond to $p < 2$ (e.g., Johnstone [21]).

4.2. *Minimaxity in sparse settings.* We now exploit Theorem 1 to prove asymptotic minimaxity of the proposed MAP estimator over various sparse balls defined above. For this purpose, we define minimax quadratic risk over a given set $\Theta_n$ in (1.1) as

$$R(\Theta_n) = \inf_{\tilde{\mu}} \sup_{\mu \in \Theta_n} E\|\tilde{\mu} - \mu\|^2$$

and examine various sparse sets $\Theta_n$, namely, $l_0$, strong and weak $l_p$-balls, where sparsity assumes that the standardized radius $\eta$ tends to zero as $n$ increases.

The general idea for establishing asymptotic minimaxity is common for all cases: for each particular setting, we find the "least favorable" sequence $\mu_0 = \mu_0(p, \eta)$ and the "equilibrium point" $k_n^* = k_n^*(p, \eta)$ that keeps balance between $\sum_{i=k_n^*+1}^n \mu_{0i}^2$ and the penalty term $P_n(k_n^*)$ on the RHS of (3.13). We show that a requirement $\pi_n(k_n^*) \geq (k_n^*/n)^{c_p k_n^*}$ for some $c_p > 0$ on the prior $\pi_n(\cdot)$ at a *single* point $k_n^*$ is sufficient for optimality of the MAP estimator $\hat{\mu}^*$. Sections 4.2.1–4.2.3 show that the MAP thresholding estimator achieves asymptotic optimality up to a constant factor in a variety of sparse settings listed above.

4.2.1. *Optimality over $l_0$-balls.* Consider an $l_0$-ball $l_0[\eta]$, where $\eta \to 0$ as $n \to \infty$. Then, by Donoho et al. [12],

$$R(l_0[\eta]) \sim \sigma^2 n\eta(2\ln \eta^{-1}),$$

where the relation "$\sim$" means that the ratio of the two sides tends to one as $n$ increases.

THEOREM 3. *Define $k_n^* = n\eta$. Let $\eta \to 0$ as $n \to \infty$ (sparsity assumption) but $n\eta \not\to 0$. If there exists a constant $c_0 > 0$ such that $\pi_n(k_n^*) \geq (k_n^*/n)^{c_0 k_n^*}$, then the MAP estimator $\hat{\mu}^*$ achieves optimality up to a constant factor, that is,*

$$\sup_{\mu \in l_0[\eta]} E\|\hat{\mu}^* - \mu\|^2 = O(n\eta(2\ln \eta^{-1})).$$

4.2.2. *Optimality over weak $l_p$-balls.* Consider a weak $l_p$-ball $m_p[\eta]$, $0 < p < 2$, and let $\eta \to 0$ as $n \to \infty$. In what follows we distinguish between sparse cases, where $n^{1/p}\eta \geq \sqrt{2\ln n}$, and *super*-sparse cases, where $n^{1/p}\eta < \sqrt{2\ln n}$. From the results of Donoho and Johnstone (e.g., Johnstone [20, 21], Donoho and Johnstone [10, 11]) it is known that

$$R(m_p[\eta]) \sim \begin{cases} \dfrac{2}{2-p}\sigma^2 n\eta^p(2\ln\eta^{-p})^{1-p/2}, & n^{1/p}\eta \geq \sqrt{2\ln n}, \\ \dfrac{2}{2-p}\sigma^2 n^{2/p}\eta^2, & n^{1/p}\eta < \sqrt{2\ln n}. \end{cases}$$

THEOREM 4. *Let $\eta \to 0$ as $n \to \infty$.*

1. *Let $n^{1/p}\eta \geq \sqrt{2\ln n}$ (sparse case). Define $k_n^* = n\eta^p(\ln\eta^{-p})^{-p/2}$. If there exists a constant $c_p > 0$ such that $\pi_n(k_n^*) \geq (k_n^*/n)^{c_p k_n^*}$, then*

$$\sup_{\mu \in m_p[\eta]} E\|\hat{\mu}^* - \mu\|^2 = O(n\eta^p(2\ln\eta^{-p})^{1-p/2}).$$

2. *Let $n^{1/p}\eta < \sqrt{2\ln n}$ (super-sparse case) but $n^{1/p}\eta \not\to 0$. If there exists a constant $c_p > 0$ such that $\pi_n(0) \geq \exp(-c_p\eta^2 n^{2/p})$, then*

$$\sup_{\mu \in m_p[\eta]} E\|\hat{\mu}^* - \mu\|^2 = O(n^{2/p}\eta^2).$$

4.2.3. *Optimality over strong $l_p$-balls.* The minimax risk over a strong $l_p$-ball, $0 < p < 2$, is the same as over the corresponding weak $l_p$-ball $m_p[\eta]$ but without the constant factor $2/(2-p)$ (Johnstone [20], Donoho and Johnstone [11]), that is,

$$R(l_p[\eta]) \sim \begin{cases} \sigma^2 n\eta^p(2\ln\eta^{-p})^{1-p/2}, & n^{1/p}\eta \geq \sqrt{2\ln n}, \\ \sigma^2 n^{2/p}\eta^2, & n^{1/p}\eta < \sqrt{2\ln n}. \end{cases}$$

THEOREM 5. *Let $\eta \to 0$ as $n \to \infty$.*

1. *Let $n^{1/p}\eta \geq \sqrt{2\ln n}$ (sparse case). Define $k_n^* = n\eta^p(\ln \eta^{-p})^{-p/2}$. If there exists a constant $c_p > 0$ such that $\pi_n(k_n^*) \geq (k_n^*/n)^{c_p k_n^*}$, then*

$$\sup_{\mu \in l_p[\eta]} E\|\hat{\mu}^* - \mu\|^2 = O(n\eta^p(2\ln \eta^{-p})^{1-p/2}).$$

2. *Let $n^{1/p}\eta < \sqrt{2\ln n}$ (super-sparse case) but $n^{1/p}\eta \not\to 0$. If there exists a constant $c_p > 0$ such that $\pi_n(0) \geq \exp(-c_p \eta^2 n^{2/p})$, then*

$$\sup_{\mu \in l_p[\eta]} E\|\hat{\mu}^* - \mu\|^2 = O(n^{2/p}\eta^2).$$

4.3. *Adaptivity.* In Sections 4.2.1–4.2.3 for sparse cases we established optimality of the MAP estimator over a given ball if the condition $\pi_n(k) \geq (k/n)^{ck}$ holds at the single "equilibrium point" $k_n^*$ depending on parameters $p$ and $\eta$ of a ball. From the results of Theorems 3–5 it follows immediately that if this condition holds for *all* $k = 1, \ldots, \kappa_n$, with some $\kappa_n < n$, the corresponding MAP estimator $\hat{\mu}^*$ is *adaptive* in the sense that it achieves optimal convergence rates *simultaneously* over an entire range of balls:

THEOREM 6 (Adaptivity). *Let $\Theta_n[\eta]$ be any of $l_0[\eta]$, $l_p[\eta]$ or $m_p[\eta]$, where $\eta \to 0$ as $n \to \infty$. If there exists $\kappa_n = o(n)$ such that $(\ln n)/\kappa_n \to 0$ as $n \to \infty$ and $\pi_n(k) \geq (k/n)^{ck}$ for all $k = 1, \ldots, \kappa_n$, and some constant $c > 0$, then, for sufficiently large $n$,*

$$\sup_{\mu \in \Theta_n[\eta]} E\|\hat{\mu}^* - \mu\|^2 = O(R(\Theta_n[\eta])) \tag{4.18}$$

*for all $0 < p < 2$ and $\eta^p \in [n^{-1}(2\ln n)^{p/2}; n^{-1}\kappa_n]$.*

*For $l_0$-balls, it is sufficient to require $\kappa_n \not\to 0$ as $n \to \infty$ [instead of $(\ln n)/\kappa_n \to 0$] in order that (4.18) hold for all $\eta \leq \kappa_n n^{-1}$ such that $n\eta \not\to 0$.*

The sufficient requirement $\pi_n(k) \geq (k/n)^{ck}$ for adaptivity established in Theorem 6 corresponds to $2k\ln(n/k)$-type penalties (see Remark 1). At the same time, the proofs of Theorems 3–5 indicate that this condition is, in fact, also "almost necessary." Thus, essentially only $2k\ln(n/k)$-type complexity penalties lead to adaptive estimation.

It is natural to find priors for which the optimality range for $\eta$ in Theorem 6 is the widest. From Theorem 6 it is clear that such priors should be of the form $\pi_n(k) \propto (k/n)^{ck}, k = 1, \ldots, \kappa_n$, where $\kappa_n = o(n)$ should be as large as possible. The function $(k/n)^k$ decreases for $k \leq n/e$ and, hence, for all $c \geq 1$, we have

$$\sum_{k=1}^{\kappa_n} \left(\frac{k}{n}\right)^{ck} \leq \sum_{k=1}^{n/e} \left(\frac{k}{n}\right)^{ck} \leq \frac{n^{1-c}}{e} \leq 1.$$

The widest possible ranges for $\eta$ in Theorem 6 are, therefore, achieved for priors of the form $\pi_n(k) = (k/n)^{ck}, k = 1, \ldots, \kappa_n$, where $c \geq 1$ and $\delta_n = \kappa_n/n$ tends to zero at an arbitrarily slow rate. The resulting ranges are $\eta \leq \delta_n$ and $\eta^p \in [n^{-1}(2\ln n)^{p/2}; \delta_n]$ for $l_0$ and $l_p$-balls, respectively, and cover the entire spectrum of sparse cases. From Lemma A.1 in the Appendix it follows that $\ln \binom{n}{k} \sim k \ln(n/k)$ for all $k = o(n)$, and the corresponding complexity penalty $P_n(k)$ in (2.10) is then

$$P_n(k) = 2\sigma^2(1 + 1/\gamma)\ln\left\{\binom{n}{k}(n/k)^{ck}(1+\gamma)^{k/2}\right\}$$

$$\sim 4\sigma^2 \tilde{c}(1 + 1/\gamma)k\left(\ln(n/k) + \frac{1}{4\tilde{c}}\ln(1+\gamma)\right),$$

where $\tilde{c} = (1/2)(c+1) \geq 1$. Such a penalty is obviously of the $2k\ln(n/k)$-type, although, by analogy, more appropriately, it should be called of the $4k\ln(n/k)$-type.

To complete this section note that, from a Bayesian viewpoint, it is also important to avoid a well-known Bayesian paradox where a prior (and, hence, a posterior) leading to an optimal estimator over a certain set has zero measure on this set. Hence, conditions on $\pi_n(k)$ should guarantee, in addition, that, with high probability, a vector $\mu$ generated according to this prior distribution falls within a given ball. We discuss this issue in more detail in the examples in Section 5.

**5. Examples.** In this section we consider three examples of $\pi_n(k)$ and establish conditions on the parameters of these distributions imposed by the general results of the previous sections.

5.1. *Binomial distribution.* Consider the binomial prior $B(n, \xi_n)$, where

$$\pi_n(k) = \binom{n}{k}\xi_n^k(1-\xi_n)^{n-k}, \qquad k = 0, \ldots, n.$$

The binomial prior suggests independent $x_i$ with $P(x_i = 1) = \xi_n, i = 1, \ldots, n$.

Assumption (A) evidently holds for any $\xi_n \leq e^{-c(\gamma)}$. We now find $\xi_n$ which satisfies the conditions of Corollary 2 and, therefore, for which the resulting MAP estimator achieves the minimal possible risk (up to a constant factor) among all available estimators in the sense of (3.17).

For $k = 0$, $\pi_n(0) = (1 - \xi_n)^n$ and in order to satisfy $\pi_n(0) \geq n^{-c}$ for some $c > 0$, $\xi_n$ should necessarily tend to zero as $n$ increases. Assumption (A) definitely holds in this case (see above). Furthermore, $(1 - \xi_n)^n = \exp\{-n\xi_n(1+o(1))\}$ and $\pi_n(0) \geq n^{-c}$ when $\xi_n \leq c_1(\ln n)/n$ for any $c_1 < c$.

On the other hand, let $\xi_n \geq n^{-c_2}$ for some $c_2 \geq 1$. Then, for all $k \geq 1$, we have

$$\pi_n(k) \geq \xi_n^k(1-\xi_n)^{n-k} \geq \xi_n^k(1-\xi_n)^n \geq n^{-c_2 k}n^{-c} \geq n^{-\tilde{c}k},$$

where $\tilde{c} = c + c_2$.

Summarizing, the validity of Corollary 2 for the binomial prior $B(n, \xi_n)$ is established for

(5.19) $$n^{-c_2} \leq \xi_n \leq c_1 \ln n/n,$$

where $c_1 > 0$ and $c_2 \geq 1$.

The condition (5.19) holds, for example, for universal thresholding, where $\xi_n \sim 1/n$, but not for the AIC criterion, where $\xi_n \sim \sqrt{\gamma}/(e + \sqrt{\gamma})$ (see the discussion in Section 2.2). Abramovich and Angelini [1] showed that for $\xi_n < \sqrt{\pi \gamma \ln n}/n$ the binomial prior leads to the Bonferroni multiple testing procedure with the familywise error rate (FWE) controlled level $\alpha_n \sim n\xi_n(\sqrt{\pi\gamma \ln(\sqrt{\gamma}/\xi_n)})^{-1} < 1$.

As we have already mentioned in Section 2.2, the binomial prior yields a *fixed* threshold $\lambda_n^2 = 2\sigma^2(1 + 1/\gamma)\ln(\frac{1-\xi_n}{\xi_n}\sqrt{1+\gamma})$ and, hence, (5.19) implies

$$\lambda_n^2 = 2\sigma^2(1+1/\gamma)(\ln \xi_n^{-1})(1+o(1)) \sim 2\sigma^2 c(\gamma)(\ln n).$$

In fact, the following proposition shows that $\xi_n$ from (5.19) also satisfies the conditions of Theorem 6 and, therefore, yields an adaptive optimal MAP estimator within the entire range of various types of sparse balls.

PROPOSITION 1. *Let $\xi_n$ satisfy (5.19). Then (4.18) holds for the resulting MAP estimator for all weak and strong $l_p$-balls with $0 < p < 2$ and $\eta^p \in [n^{-1}(2\ln n)^{p/2}; \xi_n^{c_3}]$, and for $l_0$-balls with $\eta \in [c_4 n^{-1}; \xi_n^{c_3}]$, where $0 < c_3 < 1/c_2 \leq 1$ and $c_4 > 0$ can be arbitrarily small.*

The widest possible ranges for $\eta$ covered by Proposition 1, namely, $\eta^p \in [n^{-1}(2\ln n)^{p/2}; (c_1 \ln n/n)^{c_3}]$ for $0 < p < 2$ and $\eta \in [c_4 n^{-1}; (c_1 \ln n/n)^{c_3}]$ for $p = 0$, respectively, are obtained for $\xi_n = c_1 \ln n/n$. These optimality ranges are still smaller than those for priors of the type $\pi_n(k) = (k/n)^{ck}, c \geq 1$, discussed in Section 4.3.

On the other hand, to avoid the Bayesian paradox mentioned at the end of Section 4.3, exploiting Lemma 7.1 in Abramovich et al. [4], we have for $l_0$-balls,

$$P(k > n\eta) \leq e^{-(1/4)n\xi_n \min\{|\eta/\xi_n - 1|, |\eta/\xi_n - 1|^2\}},$$

and the above probability tends to zero as $n$ increases for $\eta \geq c_5(\ln n/n)$, where $c_5 \geq 2c_1$.

For strong $l_p$-balls, define the standardized $z = \mu/\tau$ and apply Markov's inequality to get

$$P(\|\mu\|_p^p > n\eta^p) \leq e^{-n(\eta/\tau)^p} E e^{\|z\|_p^p}.$$

For the hierarchical prior model introduced in Section 2.2, $Ee^{\|z\|_p^p} = E(Ee^{\|z\|_p^p}|x) = E_{\pi_n} a_p^k$, where $a_p = Ee^{|\zeta|^p}$ and $\zeta$ is a standard normal $N(0,1)$. It is easy to verify that $e < a_p < \infty$ for $0 < p < 2$. Thus,

$$P(\|\mu\|_p^p > n\eta^p) \leq e^{-n(\eta/\tau)^p} E_{\pi_n} a_p^k. \tag{5.20}$$

For the binomial prior,

$$E_{\pi_n} a_p^k = (\xi_n a_p + 1 - \xi_n)^n = e^{n\xi_n(a_p-1)(1+o(1))}.$$

If $\eta^p \geq \tau^p \xi_n(a_p - 1) + c_6(\ln \ln n/n)$, where $c_6 > 0$ is arbitrarily small, then $P(\mu \in l_p[\eta]) \to 1$ and, therefore, $P(\mu \in m_p[\eta]) \to 1$ as well. We believe that after extra effort it is possible to somewhat relax the conditions for weak $l_p$-balls, but the resulting additional benefits are usually minor.

Combining these Bayesian admissibility results with Proposition 1, we obtain the *admissible optimality* ranges for the binomial prior $B(n, \xi_n)$ with $n^{-c_2} \leq \xi_n \leq c_1(\ln n)/n$, $c_1 > 0$ and $c_2 \geq 1$:

$$\eta \in [c_5 n^{-1} \ln n; \xi_n^{c_3}], \qquad p = 0,$$

$$\eta^p \in [\max\{n^{-1}(2\ln n)^{p/2}, \tau^p \xi_n(a_p - 1) + c_6 n^{-1} \ln \ln n\}; \xi_n^{c_3}], \qquad 0 < p < 2,$$

where $0 < c_3 < 1/c_2 \leq 1$, $c_5 \geq 2c_1$ and $c_6 > 0$ can be arbitrarily small.

5.2. *Truncated Poisson distribution.* Consider now a truncated Poisson distribution, where

$$\pi_n(k) = \frac{\lambda_n^k/k!}{\sum_{j=0}^n \lambda_n^j/j!}, \qquad k = 0, \ldots, n,$$

and $1 \leq \lambda_n \leq n$. Application of Stirling's formula and simple calculus yield the following bounds on $\pi_n(k)$:

$$\left(\frac{\lambda_n}{k}\right)^k \frac{e^{k-\lambda_n - 1/(12k)}}{\sqrt{2\pi k}} < \frac{\lambda_n^k}{k!} e^{-\lambda_n} < \pi_n(k)$$

$$< \frac{\lambda_n^k/k!}{\lambda_n^{\lambda_n}/\lambda_n!} < \left(\frac{\lambda_n}{k}\right)^{k+1/2} e^{k-\lambda_n + 1/(12\lambda_n)}. \tag{5.21}$$

From (5.21) and Lemma A.1 from the Appendix one has

$$\ln \pi_n(k) - \left(\ln \binom{n}{k} - kc(\gamma)\right)$$

$$< k \ln\left(\frac{\lambda_n e^{c(\gamma)+1}}{n}\right) - \frac{1}{2}\ln k - \lambda_n + \frac{1}{12\lambda_n} + \frac{1}{2}\ln \lambda_n. \tag{5.22}$$

The function $x - (1/12x) - (1/2)\ln x > 0$ for all $x \geq 1$ and the RHS of (5.22) is negative for all $k \geq 1$ when $\lambda_n < ne^{-(c(\gamma)+1)}$. Hence, Assumption (A) is satisfied for $\lambda_n < ne^{-(c(\gamma)+1)}$.

We now check conditions for Corollary 2. For $k=0$, one has $\pi_n(0) > e^{-\lambda_n}$ and the requirement $\pi_n(0) \geq n^{-c_1}$ of Corollary 2 is satisfied for $\lambda_n \leq c_1 \ln n$. Note that this requirement immediately yields Assumption (A).

On the other hand, let, in addition, $\lambda_n \geq n^{-c_2}$ for some $c_2 \geq 0$. For $k=1$, (5.21) implies $\pi_n(1) > \lambda_n e^{-\lambda_n} > n^{-(c_1+c_2)}$, while, for $k \geq 2$, note that $k - 1/(12k) - (1/2)\ln(2\pi k) > 0$, and, therefore, one has from (5.21)

$$\ln \pi_n(k) > k \ln\left(\frac{\lambda_n}{k}\right) - \lambda_n > k \ln\left(\frac{n^{-c_2}}{k}\right) - c_1 \ln n$$

$$> k \ln\left(\frac{n^{-c_2-c_1/k}}{n}\right) > k \ln n^{-(c_1+c_2+1)}.$$

Thus, for the truncated Poisson prior Corollary 2 holds if

(5.23) $$n^{-c_2} \leq \lambda_n \leq c_1 \ln n,$$

where $c_1 > 0$ and $c_2 \geq 0$.

In particular, Abramovich and Angelini [1] showed that for $\lambda_n < \sqrt{\pi \gamma \ln n}$ the corresponding MAP testing procedure controls the FWE at the level $\alpha_n \sim \lambda_n(\sqrt{\pi \gamma \ln(\sqrt{\gamma}n/\lambda_n)})^{-1} < 1$ and is closely related to the FWE controlling multiple testing procedures of Holm [18] and Hochberg [17].

PROPOSITION 2. *Let $\lambda_n$ satisfy (5.23). Then (4.18) holds for the resulting MAP estimator for all weak and strong $l_p$-balls with $0 < p < 2$ and $\eta^p \in [n^{-1}(2 \ln n)^{p/2}; (\lambda_n/n)^{c_3}]$, and for $l_0$-balls with $\eta \in [c_4 n^{-1}; (\lambda_n/n)^{c_3}]$, where $0 < c_3 < 1/(1+c_2)$ and $c_4 > 0$ can be arbitrarily small.*

Consider the corresponding Bayesian admissibility requirements for the truncated Poisson prior. For $l_0$-balls, in the proof of their Lemma 1, Abramovich and Angelini [1] showed that, with $\lambda_n = o(n)$ for any $\delta_n = o(n)$,

$$P(k \geq \lambda_n + \delta_n) \leq Cnu_n,$$

where $u_n = e^{\delta_n}/(1+\delta_n/\lambda_n)^{\lambda_n+\delta_n+1/2}$ and, therefore, $\ln u_n < \delta_n(1-\ln(\delta_n/\lambda_n))$. In particular, set $\delta_n = \max(\ln n, e^\zeta \lambda_n)$, where $\zeta > 2$. Then

$$P(k \geq \lambda_n + \delta_n) \leq Cne^{-\delta_n(\zeta-1)} \leq Cn^{-(\zeta-2)} \to 0$$

and, hence, $P(\mu \in l_0[\eta]) \to 1$ for $\eta \geq (\lambda_n + \delta_n)/n$. For $\lambda_n$ satisfying (5.23), $P(\mu \in l_0[\eta]) \to 1$ holds for $\eta \geq c_5(\ln n/n)$, where $c_5 > 2\max(1, e^2 c_1)$.

For $l_p$-balls, exploiting (5.20) for the truncated Poisson prior and applying Stirling's formula, one derives

$$E_{\pi_n} a_p^k = \frac{\sum_{k=0}^n a_p^k \lambda_n^k/k!}{\sum_{j=0}^n \lambda_n^j/j!} \leq \frac{\sum_{k=0}^\infty a_p^k \lambda_n^k/k!}{\lambda_n^{\lambda_n}/\lambda_n!} \leq e^{\lambda_n(a_p-1)}\sqrt{2\pi\lambda_n}$$

and, therefore, for $\eta^p \geq \tau^p(\lambda_n/n)(a_p - 1) + c_6(\ln\ln n/n)$, where $c_6 > 0$ is arbitrarily small, both $P(\mu \in l_p[\eta])$ and $P(\mu \in m_p[\eta])$ tend to one.

The resulting admissible optimality ranges for $\eta$ for the truncated Poisson prior with $n^{-c_2} \leq \lambda_n \leq c_1 \ln n$, $c_1 > 0$ and $c_2 \geq 0$ are then given by

$$\eta \in [c_5(\ln n/n); (\lambda_n/n)^{c_3}], \qquad p = 0,$$

$$\eta^p \in [n^{-1}\max\{(2\ln n)^{p/2}, \tau^p \lambda_n(a_p - 1) + c_6 \ln\ln n\}; (\lambda_n/n)^{c_3}], \qquad 0 < p < 2,$$

where $0 < c_3 < 1/(1 + c_2)$, $c_5 > 2\max(1, e^2 c_1)$, and $c_6 > 0$ can be arbitrarily small.

Strong similarity between the results for truncated Poisson and binomial priors with $\xi_n = \lambda_n/n$ is not surprising and is due to the well-known asymptotic relations between Poisson and binomial distributions.

5.3. *Reflected truncated Poisson distribution.* Finally, consider briefly a "reflected" truncated Poisson distribution

$$(5.24) \qquad \pi_n(k) = \frac{(n - \lambda_n)^{n-k}/(n-k)!}{\sum_{j=0}^n (n - \lambda_n)^{n-j}/(n-j)!}, \qquad k = 0, \ldots, n,$$

and let $\lambda_n = o(n)$ but $\lambda_n/\sqrt{n \ln n} \to \infty$ as $n$ increases. The motivation for such a type of prior and specific choice of $\lambda_n$ comes from the fact that the corresponding MAP testing procedure mimics the FDR controlling procedures of Benjamini and Hochberg [7] and Sarkar [23] (Abramovich and Angelini [1]). In particular, Abramovich and Angelini ([1], Lemma 2) showed that, almost surely, $k = \lambda_n(1 + o(1))$ or, more precisely, $|k - \lambda_n| \leq \sqrt{c_7 n \ln n}$, where $c_7 > 4$.

The adaptivity results of Theorem 6 are somewhat irrelevant for such a narrow range of possible $k$ since the "equilibrium point" $k_n^* = \lambda_n(1 + o(1)) = o(n)$ in Theorems 3–5 becomes essentially known. To apply Theorems 3–5, we need the following lemma.

LEMMA 1. *Consider the reflected truncated Poisson prior $\pi_n(k)$ with $\lambda_n = o(n)$, $\lambda_n/\sqrt{n \ln n} \to \infty$ as $n \to 0$. For $k = \lambda_n(1 + o(1))$, there exists $c > 0$ such that $\pi_n(k) \geq (k/n)^{ck}$.*

Based on the results of Lemma 1, we can identify the radius $\eta_0$ of the balls, where the resulting MAP estimator $\hat{\mu}^*$ is optimal. For $l_0$-balls, Theorem 3 yields $\eta_0 = (\lambda_n/n)(1 + o(1))$. Similarly, for $0 < p < 2$, applying Theorem 4 and Theorem 5, we obtain the result that the corresponding $\eta_0$ satisfies $\eta_0^p(\ln \eta_0^{-p})^{-p/2} = (\lambda_n/n)(1 + o(1))$.

We now show that there is no Bayesian paradox in this case and a vector $\mu$ generated from $\pi_n(k)$ falls with high probability within the corresponding balls of radius $\eta_0$. For $l_0$-balls, it follows since, almost surely, $k = \lambda_n(1 +$ $o(1)) \sim n\eta_0$. For $l_p$-balls, $0 < p < 2$, Abramovich and Angelini [1] proved that for the reflected truncated Poisson prior (5.24) we have $Ek = \lambda_n(1 + o(1))$. Then, Markov's inequality implies

$$P(\|\mu\|_p^p > n\eta_0^p) \leq \frac{E\|\mu\|_p^p}{n\eta_0^p} = \frac{\tau^p \nu_p E_{\pi_n} k}{n\eta_0^p} = \tau^p \nu_p (\ln \eta_0^{-p})^{-p/2} \to 0,$$

where $\nu_p$ is the $p$th absolute moment of the standard normal distribution.

TABLE 1
*AMSE of various thresholding estimators*

| $\xi\%$ | 0.5% | | | 5% | | | 50% | | |
|---|---|---|---|---|---|---|---|---|---|
| $\tau$ | 3 | 5 | 7 | 3 | 5 | 7 | 3 | 5 | 7 |
| Bin | 0.0192 | 0.0194 | 0.0172 | 0.1496 | 0.1372 | 0.1245 | 0.8929 | 0.8196 | 0.7796 |
| Pois1 | 0.0192 | 0.0194 | 0.0176 | 0.1497 | 0.1374 | 0.1246 | 0.9157 | 0.8245 | 0.7780 |
| Pois2 | 0.0187 | 0.0195 | 0.0173 | 0.1564 | 0.1389 | 0.1256 | 0.9687 | 0.8271 | 0.7795 |
| EbayesThresh | 0.0194 | 0.0189 | 0.0181 | 0.1556 | 0.1379 | 0.1248 | 1.8942 | 1.3473 | 0.8450 |
| universal | 0.1012 | 0.0998 | 0.0990 | 0.1694 | 0.1600 | 0.1495 | 1.9447 | 2.5474 | 2.7720 |

**6. Some simulation results.** A short simulation study was carried out to investigate the performance of several MAP estimators.

The data was generated according to the model (1.1) with the sample size $n = 1000$. In $\xi\%$ percent of cases $\mu_i$ were randomly sampled from $N(0, \tau^2)$, and otherwise $\mu_i = 0$. The parameter $\xi$ controls the sparsity of the true signal $\mu$, while $\tau$ reflects its energy. We considered $\xi = 0.5\%$, 5% and 50% corresponding, respectively, to super-sparse, sparse and dense cases, and $\tau = 3, 5, 7$. For each combination of values of $\xi$ and $\tau$, the number of replications was 100. The true values of $\sigma$, $\tau$ and $\xi$ were assumed *unknown* in simulations and were estimated from the data by the EM-algorithm of Abramovich and Angelini [1]. Our simulation study also confirmed the efficiency of their parameter estimation procedure.

We tried three MAP estimators corresponding to the priors considered in Section 5, namely, binomial $B(n, \xi)$ (Bin), truncated Poisson (Pois1) and reflected truncated Poisson (Pois2) with $\lambda = n\xi$. In addition, we compared performances of the above listed MAP estimators with the universal thresholding of Donoho and Johnstone [9] and the hard thresholding EbayesThresh estimator of Johnstone and Silverman [22] with a Cauchy prior.

Table 1 summarizes mean squared errors averaged over 100 replications (AMSE) of various methods. Standard errors in all cases were of order several percent of the corresponding AMSE.

The performance of all methods naturally improves as $\tau$ increases. As is typical for any thresholding procedure, all of them are less efficient for dense cases. Nonadaptive universal thresholding consistently yields the worst results. All MAP estimators behave similarly, indicating the robustness of the MAP testimation to the choice of the prior $\pi_n(\cdot)$. They are comparable with EbayesThresh for very sparse and sparse cases, but strongly outperform the latter for dense signals. Partially this is explained by the poor behavior of the MAD estimate of $\sigma$ used in this case by EbayesThresh. However, even after substituting the true $\sigma$ in EbayesThresh, MAP estimators still remained preferable.

**7. Concluding remarks.** In this paper we have considered a Bayesian approach to a high-dimensional normal means problem. The proposed hierarchical prior is based on assuming a prior distribution $\pi_n(\cdot)$ on the number of nonzero entries of the unknown means vector. The resulting Bayesian MAP "testimator" leads to a hard thresholding rule and, from a frequentist viewpoint corresponds to penalized likelihood estimation with a complexity penalty depending on $\pi_n(\cdot)$. Specific choices of $\pi_n(\cdot)$ lead to several well-known complexity penalties. In particular, we have discussed the relationship between MAP testimation and $2k\ln(n/k)$-type penalization recently considered in a series of papers. We have investigated the optimality of MAP estimators and established their adaptive minimaxity in various sparse settings.

In practice, the unknown parameters of the prior and the noise variance can be efficiently estimated by the EM algorithm. The simulation study presented illustrates the theoretical results and shows the robustness properties of the MAP testimation procedure to the choice of $\pi_n(\cdot)$.

We believe that the proposed Bayesian approach for recovering a high-dimensional vector from white noise can be extended to various non-Gaussian settings and model selection problems, although appropriate adjustments are needed for each specific problem at hand.

## APPENDIX

**A.1. Proof of Theorem 1.** We show that, for a prior satisfying Assumption (A), the corresponding penalty $P_n(k)$ in (2.10) belongs to the general class of penalties considered in Birgé and Massart [8]. In particular, we verify that it satisfies conditions (3.3) and (3.4) of their Theorem 2 and then use it directly to obtain the upper bound in (3.13).

In our notation the conditions (3.3) and (3.4) of Birgé and Massart [8] correspond respectively to

$$\text{(A.1)} \qquad \sum_{k=1}^{n} \binom{n}{k} e^{-kL_{k,n}} < \infty$$

and

(A.2)    $(1+1/\gamma)(2L_{k,n} + \ln(1+\gamma)) \geq c(1+\sqrt{2L_{k,n}})^2$,     $k=1,\ldots,n$,

for some $c > 1$, where the weights $L_{k,n}$ were defined in (3.15). In fact, Birgé and Massart [8] require their (3.4) for $k = 0$ as well. However, note that $P_n(0) = 2\sigma^2(1+1/\gamma)\ln \pi_n(0)^{-1} \geq 0$ and, hence, this condition always holds for $k = 0$.

The condition (A.1) follows immediately from the definition of $L_{k,n}$,

$$\sum_{k=1}^n \binom{n}{k} e^{-kL_{k,n}} = 1 - \pi_n(0) < \infty.$$

We now turn to (A.2). Consider $k \geq 1$. Let $t = \sqrt{L_{k,n}}$. The condition (A.2) is then equivalent to the quadratic inequality

(A.3)     $2(1+1/\gamma - c)t^2 - 2\sqrt{2}ct + (1+1/\gamma)\ln(1+\gamma) - c \geq 0$.

We now find $c > 1$, for which (A.3) holds for all $t$ such that the corresponding $L_{k,n}$ satisfy Assumption (A). For the determinant $\Delta$ of (A.3), one has

$$\frac{\Delta}{4} = 2c^2 - 2(1+1/\gamma - c)((1+1/\gamma)\ln(1+\gamma) - c)$$

$$= 2(1+1/\gamma)(c(\ln(1+\gamma)+1) - (1+1/\gamma)\ln(1+\gamma)).$$

Note that $\ln(1+\gamma) \leq \gamma$ and, therefore, $\ln(1+\gamma)(1+1/\gamma) \leq \ln(1+\gamma) + 1$. Hence, $\Delta \geq 0$ for any $c > 1$. If, in addition, $c < 1 + 1/\gamma$, then (A.3) holds for all $t \geq t_*$, where $t_*$ is the largest root of the quadratic polynomial on the left-hand side of (A.3),

(A.4)
$$t_* = \frac{c + \sqrt{1+1/\gamma}\sqrt{c(\ln(1+\gamma)+1) - (1+1/\gamma)\ln(1+\gamma)}}{\sqrt{2}(1+1/\gamma - c)}$$

$$< \frac{c + 1 + 1/\gamma}{\sqrt{2}(1+1/\gamma - c)}.$$

Setting $c = 1 + 1/(2\gamma)$, from (A.4) one has $t_* < 2\sqrt{2}(\gamma + 3/4)$. On the other hand, Assumption (A) implies $L_{k,n} \geq c(\gamma) = 8(\gamma + 3/4)^2$ and, therefore, $t = \sqrt{L_{k,n}} \geq 2\sqrt{2}(\gamma + 3/4) > t_*$. Thus, $c = 1 + 1/(2\gamma)$ guarantees the condition (A.3) and the equivalent original condition (A.2).

**A.2. Proof of Theorem 2.** Consider first $k \geq 1$. For this case Assumption (A) implies $L_{k,n} \geq c(\gamma) \geq c(0) > 1$. From Theorem 1, we then have

$$\rho(\hat{\mu}^*, \mu) \leq c_0(\gamma)(1+1/\gamma)$$

$$\times \inf_{1 \leq k \leq n}\left\{\sum_{i=k+1}^n \mu_{(i)}^2 + \sigma^2 k(2L_{k,n} + \ln(1+\gamma))\right\} + c_1(\gamma)\sigma^2$$

$$
\begin{aligned}
&\leq c_0(\gamma)(1+1/\gamma)(2L_n^* + \ln(1+\gamma)) \\
&\quad \times \inf_{1\leq k\leq n}\left\{\sum_{i=k+1}^n \mu_{(i)}^2 + k\sigma^2\right\} + c_1(\gamma)\sigma^2 \\
&\leq c_2(\gamma)(2L_n^* + \ln(1+\gamma))\left\{\inf_{1\leq k\leq n}\left(\sum_{i=k+1}^n \mu_{(i)}^2 + k\sigma^2\right) + \sigma^2\right\}.
\end{aligned}
\tag{A.5}
$$

On the other hand, Theorem 1 implies

$$
\begin{aligned}
\rho(\hat{\mu}^*,\mu) &\leq c_0(\gamma)\left\{\sum_{i=1}^n \mu_i^2 + 2\sigma^2(1+1/\gamma)\ln\pi_n^{-1}(0)\right\} + c_1(\gamma)\sigma^2 \\
&\leq c_0(\gamma)\left\{\sum_{i=1}^n \mu_i^2 + \sigma^2(1+1/\gamma)(L_{0,n} + 0.5\ln(1+\gamma))\right\} + c_1(\gamma)\sigma^2.
\end{aligned}
$$

Define $\tilde{c}_0(\gamma) = c_0(\gamma)/\ln(1+\gamma)$ and $\tilde{c}_1(\gamma) = 2c_1(\gamma)/\ln(1+\gamma)$. Obviously,

$$\tilde{c}_0(\gamma)(2L_{0,n} + \ln(1+\gamma)) > c_0(\gamma)$$

and

$$\tilde{c}_1(\gamma)(2L_{0,n} + \ln(1+\gamma)) > 2c_1(\gamma).$$

Hence,

$$
\begin{aligned}
\rho(\hat{\mu}^*,\mu) &\leq \tilde{c}_0(\gamma)(2L_{0,n}+\ln(1+\gamma))\sum_{i=1}^n \mu_i^2 \\
&\quad + c_0(\gamma)(1+1/\gamma)(2L_{0,n}+\ln(1+\gamma))\sigma^2/2 \\
&\quad + \tilde{c}_1(\gamma)(2L_{0,n}+\ln(1+\gamma))\sigma^2/2 \\
&\leq c_2(\gamma)(2L_n^*+\ln(1+\gamma))\left\{\sum_{i=1}^n \mu_i^2 + \sigma^2\right\}.
\end{aligned}
\tag{A.6}
$$

Combining (A.5) and (A.6), we have

$$
\begin{aligned}
\rho(\hat{\mu}^*,\mu) &\leq c_2(\gamma)(2L_n^*+\ln(1+\gamma))\left\{\inf_{0\leq k\leq n}\left(\sum_{i=k+1}^n \mu_{(i)}^2 + k\sigma^2\right) + \sigma^2\right\} \\
&= c_2(\gamma)(2L_n^*+\ln(1+\gamma))\left\{\sum_{i=1}^n \min(\mu_i^2,\sigma^2) + \sigma^2\right\},
\end{aligned}
$$

which completes the proof.

**A.3. Proof of Theorem 3.** We start with Lemma A.1, which will be used throughout the following proofs.

LEMMA A.1.

1. $\ln \binom{n}{k} \geq k \ln(n/k), k = 1, \ldots, n$.
2. Let $n/k \to \infty$ as $n \to \infty$. Then for any constant $c > 1$ for sufficiently large $n$, $\ln \binom{n}{k} \leq ck \ln(n/k)$.

PROOF. The first statement of Lemma A.1 follows immediately from the trivial inequality

$$\binom{n}{k} = \prod_{j=0}^{k-1} \frac{n-j}{k-j} \geq \left(\frac{n}{k}\right)^k.$$

To prove the second statement, note that, using Stirling's formula, one has

(A.7)
$$\binom{n}{k} \leq \left(\frac{n}{e}\right)^n \left(\frac{e}{n-k}\right)^{n-k} \left(\frac{e}{k}\right)^k$$
$$= \left(\frac{n}{k}\right)^k \left(\frac{n}{n-k}\right)^{n-k} < \left(\frac{n}{k}\right)^k \left(\frac{n}{n-k}\right)^n.$$

Since $(1 - k/n)^{-n/k} \to e$ as $n/k \to \infty$, for any $c > 1$ for sufficiently large $n$,

$$\left(\frac{n}{n-k}\right)^n = ((1-k/n)^{-n/k})^k < \left(\frac{n}{k}\right)^{(c-1)k}.$$

Thus, from (A.7), $\binom{n}{k} < (n/k)^{ck}$ for sufficiently large $n$. □

Now we return to the proof of Theorem 3. Evidently, for any $\mu \in l_0[\eta]$, $\mu_{(i)} = 0, i > k^* = n\eta$. Since $k^* = o(n)$, from the general upper bound for the risk established in Corollary 1, it follows that

$E\|\hat{\mu}^* - \mu\|^2$
$$\leq c_0(\gamma) 2\sigma^2 (1 + 1/\gamma) \left( \ln\left\{ \binom{n}{n\eta} \pi_n^{-1}(n\eta) \right\} + \frac{n\eta}{2} \ln(1+\gamma) \right) + c_1(\gamma)\sigma^2.$$

From Lemma A.1,

$$\ln\left\{ \binom{n}{n\eta} \pi_n^{-1}(n\eta) \right\} \geq \ln \binom{n}{n\eta} \geq n\eta \ln \eta^{-1} \gg n\eta \ln(1+\gamma)$$

when $\eta \to 0$ as $n \to \infty$. On the other hand, under the conditions of Theorem 3, Lemma A.1 implies

$$\ln\left\{ \binom{n}{n\eta} \pi_n^{-1}(n\eta) \right\} \leq \tilde{c} n\eta \ln \eta^{-1}$$

for sufficiently large $n$. Summarizing, one has $E\|\hat{\mu}^* - \mu\|^2 \leq \tilde{c}_2(\gamma)\sigma^2 n\eta \ln \eta^{-1}$.

**A.4. Proof of Theorem 4.** Define a "least-favorable" sequence $\mu_{0i} = \eta(n/i)^{1/p}$, $i = 1, \ldots, n$, that maximizes $\sum_{i=k+1}^{n} \mu_i^2$ over $\mu \in m_p[\eta]$ for any $k = 0, \ldots, n-1$. For $k \geq 1$,

$$\text{(A.8)} \quad \sum_{i=k+1}^{n} \mu_{0i}^2 \leq \eta^2 n^{2/p} \int_k^\infty x^{-2/p}\,dx = \frac{p}{2-p}\eta^2 n^{2/p} k^{1-2/p},$$

while, for $k = 0$,

$$\sum_{i=1}^{n} \mu_{0i}^2 \leq \eta^2 n^{2/p} \zeta(2/p),$$

where $\zeta(\cdot) < \infty$ is the Riemann Zeta-function.

1. $n^{1/p}\eta \geq \sqrt{2\ln n}$ (sparse case).

In this case, $1 < k_n^* = o(n)$ and from Corollary 1, Lemma A.1 and (A.8), one has

$$E\|\hat{\mu}^* - \mu\|^2 \leq c_0(\gamma)\Bigg\{\sum_{i=k_n^*+1}^{n} \mu_{0i}^2 + 2\sigma^2(1+1/\gamma)$$

$$\times \left(\ln\binom{n}{k_n^*} + \ln \pi_n^{-1}(k_n^*) + \frac{k_n^*}{2}\ln(1+\gamma)\right)\Bigg\}$$

$$+ c_1(\gamma)\sigma^2$$

$$\leq \tilde{c}_0(\gamma)\Bigg\{\frac{p}{2-p}\eta^2 n^{2/p}(k_n^*)^{1-2/p}$$

$$+ \tilde{c}_1(\gamma)\sigma^2(k_n^* \ln(n/k_n^*) + \ln \pi_n^{-1}(k_n^*))\Bigg\}.$$

To complete the proof for this case, note that $\eta^2 n^{2/p}(k_n^*)^{1-2/p} = n\eta^p(\ln \eta^{-p})^{1-p/2}$ and under the conditions on $\pi_n(k_n^*)$ of Theorem 4,

$$k_n^* \ln(n/k_n^*) + \ln \pi_n^{-1}(k_n^*) \leq (c_p + 1)k_n^* \ln(\eta^{-p}(\ln \eta^{-p})^{p/2})$$

$$= O(n\eta^p(\ln \eta^{-p})^{1-p/2}).$$

2. $n^{1/p}\eta < \sqrt{2\ln n}$ (super-sparse case).

In this case Corollary 1 and conditions on $\pi_n(0)$ imply

$$E\|\hat{\mu}^* - \mu\|^2 \leq c_0(\gamma)\Bigg\{\sum_{i=1}^{n}\mu_{0i}^2 + 2\sigma^2(1+1/\gamma)\ln \pi_n^{-1}(0)\Bigg\} + c_1(\gamma)\sigma^2$$

$$\leq c_0(\gamma)\{\eta^2 n^{2/p}\zeta(2/p) + 2\sigma^2(1+1/\gamma)c_p\eta^2 n^{2/p}\} + c_1(\gamma)\sigma^2$$

$$= O(\eta^2 n^{2/p}).$$

**A.5. Proof of Theorem 5.** First, we find the "least favorable" sequence $\mu_0$ that maximizes $\sum_{i=k+1}^{n} \mu_{(i)}^2$ over $\mu \in l_p[\eta]$ for a given $k = 0, \ldots, n-1$. Applying Lagrange multipliers, after some algebra one has

$$\sum_{i=k+1}^{n} \mu_{0(i)}^2 \leq \left(\frac{2-p}{2}\right)^{2/p} \frac{p}{2-p} \eta^2 n^{2/p} k^{1-2/p}$$

for $k \geq 1$ and

$$\sum_{i=1}^{n} \mu_{0(i)}^2 \leq \eta^2 n^{2/p}$$

for $k = 0$. The rest of the proof therefore repeats the proof of Theorem 4.

**A.6. Proof of Proposition 1.** First, note that $\pi_n(1) = n\xi_n(1-\xi_n)^{n-1}$ and the condition $\pi_n(1) \geq n^{-c}$ in Theorem 6 requires that $\xi_n \to 0$ as $n \to \infty$. In particular, it implies

$$(1-\xi_n)^n = ((1-\xi_n)^{-1/\xi_n})^{-n\xi_n} = \exp\{-n\xi_n(1+o(1))\}$$

and, using Lemma A.1, we then have

$$\pi_n(k) \geq \left(\frac{n}{k}\right)^k \xi_n^k (1-\xi_n)^n = \left(\frac{n\xi_n}{k} \exp\left\{-\frac{n\xi_n}{k}(1+o(1))\right\}\right)^k.$$

To satisfy $\pi_n(k) \geq (k/n)^{ck}$, it is sufficient to have

(A.9) $$\frac{n\xi_n}{k}(1+o(1)) - \ln\left(\frac{n\xi_n}{k}\right) \leq c\ln\left(\frac{n}{k}\right).$$

Recall that $n^{-c_2} \leq \xi_n \leq c_1(\ln n)/n$, where $c_1 > 0$ and $c_2 \geq 1$. Define $\kappa_n = n\xi_n^{c_3}$, where $0 < c_3 < 1/c_2$. Obviously, $\kappa_n/n \to 0$, $\kappa_n > n^{1-c_2 c_3} = n^\delta, \delta > 0$, and, therefore, $(\ln n)/\kappa_n \to 0$ as $n \to \infty$.

For $1 \leq k \leq n\xi_n$, using the monotonicity of the function $-x \ln x$ for $x \leq 1/e$, we have

$$\frac{k}{n}\ln\left(\frac{n}{k}\right) \geq \frac{\ln n}{n} \geq c_1^{-1}\xi_n$$

and, therefore,

$$\frac{n\xi_n}{k}(1+o(1)) - \ln\left(\frac{n\xi_n}{k}\right) \leq \frac{n\xi_n}{k}(1+o(1)) < c_1(1+o(1))\ln\left(\frac{n}{k}\right)$$

which yields (A.9).

On the other hand, for all $n\xi_n < k \leq \kappa_n$, we have $\xi_n \geq (k/n)^{1/c_3} = (k/n)^{1+\tilde{c}}$, where $\tilde{c} > 0$, which yields $(n\xi_n/k) \geq (k/n)^{\tilde{c}}$. Thus,

$$\frac{n\xi_n}{k}(1+o(1)) - \ln\left(\frac{n\xi_n}{k}\right) < (1+o(1)) + \tilde{c}\ln\left(\frac{n}{k}\right)$$

and (A.9) holds.

Applying Theorem 6 for $\kappa_n = n\xi_n^{c_3}$ completes the proof.

**A.7. Proof of Proposition 2.** Define $\kappa_n = n(\lambda_n/n)^{c_3}$, where $0 < c_3 < 1/(1+c_2)$. Obviously, $\kappa_n > n^{1-(1+c_2)c_3} = n^\delta$, where $0 < \delta < 1$ and, therefore, $(\ln n)/\kappa_n \to 0$ and $\kappa_n/n = (\lambda_n/n)^{c_3} < (c_1 \ln n/n)^{c_3} \to 0$ as $n \to \infty$.

For $k = 1$, we have shown in Section 5.2 that $\pi_n(1) \geq n^{-(c_1+c_2)}$. For $k \geq 2$, exploit positivity of the function $k - 1/(12k) - (1/2)\ln(2\pi k)$ to obtain from (5.21)

$$\ln \pi_n(k) > k \ln\left(\frac{\lambda_n}{k}\right) - \lambda_n.$$

The rest of the proof essentially repeats the proof of Proposition 1 starting from (A.9) with $\lambda_n = n\xi_n$ and without $o(1)$.

**A.8. Proof of Lemma 1.** Applying Stirling's formula for large $\lambda_n$ and $k = \lambda_n(1 + o(1))$, after simple calculation, one has

$$\begin{aligned}
\ln \pi_n(k) &> (n-k)\ln\frac{n-\lambda_n}{n-k} + \lambda_n - k - \frac{1}{12(n-k)} - \frac{1}{2}\ln(2\pi(n-k)) \\
&= o(\lambda_n)\ln\left(\left(1 + \frac{o(\lambda_n)}{n-\lambda_n-o(\lambda_n)}\right)^{(n-\lambda_n-o(\lambda_n))/o(\lambda_n)}\right) \\
&\quad - o(\lambda_n) - \frac{1}{2}\ln(n - \lambda_n - o(\lambda_n)) \\
&= o(\lambda_n) - \frac{1}{2}\ln(n-\lambda_n) - \frac{1}{2}\ln\left(\frac{n-\lambda_n-o(\lambda_n)}{n-\lambda_n}\right) \\
&= o(\lambda_n) - \frac{1}{2}\ln(n-\lambda_n).
\end{aligned}$$

On the other hand, $ck\ln(k/n) = c\lambda_n(1+o(1))\ln(\lambda_n/n) + o(\lambda_n)$. Thus, to prove Lemma 1, it is sufficient to show that

$$(A.10) \qquad \tfrac{1}{2}\ln(n-\lambda_n) \leq \tilde{c}\lambda_n \ln(n/\lambda_n)$$

for some $\tilde{c} > 0$.

Denote $g_1(\lambda_n) = \frac{1}{2}\ln(n-\lambda_n)$ and $g_2(\lambda_n) = \lambda_n \ln(n/\lambda_n)$. Note that $g_1(\lambda_n)$ decreases while $g_2(\lambda_n)$ increases for $\lambda_n < n/e$, and $g_1(1) < g_2(1)$. Then, for any $\tilde{c} \geq 1$, one has $\tilde{c}g_2(\lambda_n) > \tilde{c}g_2(1) \geq g_1(1) > g_1(\lambda_n)$, which proves (A.10).

**Acknowledgments.** The authors would like to thank Claudia Angelini for the MATLAB codes and Nasia Petsa for pointing out several mistakes in the earlier versions of the proofs. Helpful comments of the Editor Jianqing Fan, the Associate Editor and two anonymous referees, which helped to improve the paper substantially, are gratefully acknowledged.

# REFERENCES

[1] ABRAMOVICH, F. and ANGELINI, C. (2006). Bayesian maximum a posteriori multiple testing procedure. *Sankhyā* **68** 436–460. MR2322194

[2] ABRAMOVICH, F. and BENJAMINI, Y. (1995). Thresholding of wavelet coefficients as a multiple hypotheses testing procedure. In *Wavelets and Statistics. Lecture Notes in Statist.* **103** 5–14. Springer, New York.

[3] ABRAMOVICH, F. and BENJAMINI, Y. (1996). Adaptive thresholding of wavelet coefficients. *Comput. Statist. Data Anal.* **22** 351–361. MR1411575

[4] ABRAMOVICH, F., BENJAMINI, Y., DONOHO, D. L. and JOHNSTONE, I. M. (2006). Adapting to unknown sparsity by controlling the false discovery rate. *Ann. Statist.* **34** 584–653. MR2281879

[5] AKAIKE, H. (1973). Information theory and an extension of the maximum likelihood principle. In *Second International Symposium on Information Theory* (B. N. Petrov and F. Csáki, eds.) 267–281. Akadémiai Kiadó, Budapest. MR0483125

[6] ANTONIADIS, A. and FAN, J. (2001). Regularization of wavelet approximations (with discussion). *J. Amer. Statist. Assoc.* **96** 939–967. MR1946364

[7] BENJAMINI, Y. and HOCHBERG, Y. (1995). Controlling the false discovery rate: A practical and powerful approach to multiple testing. *J. Roy. Statist. Soc. Ser. B* **57** 289–300. MR1325392

[8] BIRGÉ, L. and MASSART, P. (2001). Gaussian model selection. *J. Eur. Math. Soc.* **3** 203–268. MR1848946

[9] DONOHO, D. L. and JOHNSTONE, I. M. (1994). Ideal spatial adaptation via wavelet shrinkage. *Biometrika* **81** 425–455. MR1311089

[10] DONOHO, D. L. and JOHNSTONE, I. M. (1994). Minimax risk over $\ell_p$-balls for $\ell_q$-error. *Probab. Theory Related Fields* **99** 277–303. MR1278886

[11] DONOHO, D. L. and JOHNSTONE, I. M. (1996). Neo-classical minimax problems, thresholding and adaptive function estimation. *Bernoulli* **2** 39–62. MR1394051

[12] DONOHO, D. L., JOHNSTONE, I. M., HOCH, J. C. and STERN, A. S. (1992). Maximum entropy and the nearly black object (with discussion). *J. Roy. Statist. Soc. Ser. B* **54** 41–81. MR1157714

[13] FAN, J. and LI, R. (2001). Variable selection via nonconcave penalized likelihood and its oracle properties. *J. Amer. Statist. Assoc.* **96** 1348–1360. MR1946581

[14] FOSTER, D. and GEORGE, E. (1994). The risk inflation criterion for multiple regression. *Ann. Statist.* **22** 1947–1975. MR1329177

[15] FOSTER, D. and STINE, R. (1999). Local asymptotic coding and the minimum description length. *IEEE Trans. Inform. Theory* **45** 1289–1293. MR1686271

[16] FRANK, I. E. and FRIEDMAN, J. H. (1993). A statistical view of some chemometrics regression tools (with discussion). *Technometrics* **35** 109–148.

[17] HOCHBERG, Y. (1988). A sharper Bonferroni procedure for multiple tests of significance. *Biometrika* **75** 800–802. MR0995126

[18] HOLM, S. (1979). A simple sequentially rejective multiple test procedure. *Scand. J. Statist.* **6** 65–70. MR0538597

[19] HUNTER, D. R. and LI, R. (2005). Variable selection using MM algorithms. *Ann. Statist.* **33** 1617–1642. MR2166557

[20] JOHNSTONE, I. M. (1994). Minimax Bayes, asymptotic minimax and sparse wavelet priors. In *Statistical Decision Theory and Related Topics V* (S. Gupta and J. Berger, eds.) 303–326. Springer, New York. MR1286310

[21] JOHNSTONE, I. M. (2002). Function estimation and Gaussian sequence models. Unpublished manuscript.

[22] JOHNSTONE, I. M. and SILVERMAN, B. W. (2004). Needles and straw in haystacks: Empirical Bayes estimates of possibly sparse sequences. *Ann. Statist.* **32** 1594–1649. MR2089135
[23] SARKAR, S. K. (2002). Some results on false discovery rate in stepwise multiple testing procedures. *Ann. Statist.* **30** 239–257. MR1892663
[24] SCHWARZ, G. (1978). Estimating the dimension of a model. *Ann. Statist.* **6** 461–464. MR0468014
[25] TIBSHIRANI, R. (1996). Regression shrinkage and selection via the lasso. *J. Roy. Statist. Soc. Ser. B* **58** 267–288. MR1379242
[26] TIBSHIRANI, R. and KNIGHT, K. (1999). The covariance inflation criterion for adaptive model selection. *J. R. Stat. Soc. Ser. B Stat. Methodol.* **61** 529–546. MR1707859

F. ABRAMOVICH
DEPARTMENT OF STATISTICS
  AND OPERATIONS RESEARCH
TEL AVIV UNIVERSITY
TEL AVIV 69978
ISRAEL
E-MAIL: felix@post.tau.ac.il

V. GRINSHTEIN
DEPARTMENT OF MATHEMATICS
THE OPEN UNIVERSITY OF ISRAEL
RAANANA 43107
ISRAEL
E-MAIL: vadimg@oumail.openu.ac.il

M. PENSKY
DEPARTMENT OF MATHEMATICS
UNIVERSITY OF CENTRAL FLORIDA
ORLANDO, FLORIDA 32816-1353
USA
E-MAIL: mpensky@pegasus.cc.ucf.edu